\documentclass[12pt,reqno]{amsart} 
\usepackage[latin1]{inputenc}
\usepackage{amsmath, amsthm, amssymb,color,graphicx,bm,psfrag,float,
  verbatim, tikz-cd, fullpage}
\usepackage[modulo, right]{lineno}
\numberwithin{equation}{section}

\usepackage[all]{xy} 
\usepackage{fullpage}
\usepackage{hyperref} 
\newcommand{\ra}{\rightarrow} 
\newcommand{\lra}{\longrightarrow} 
\newcommand{\conic}{\mathcal K}
\newcommand{\complex}{\mathbf C}
\newcommand{\textC}{\mathcal{C}} 
\newcommand{\proj}{\mathbb P} 
\newcommand{\ltr}{\text{\sc ltr}}
\newcommand{\bA}{\mathbb A} 
\newcommand{\bB}{\mathbb B} 
\newcommand{\bC}{\mathbb C} 
\newcommand{\bD}{\mathbb D} 
\newcommand{\bE}{\mathbb E} 
\newcommand{\bF}{\mathbb F}
\newcommand{\bT}{\mathbb T}
\newcommand{\bV}{\mathbb V}
\newcommand{\bbX}{\mathbb X}
\newcommand{\fraka}{\mathbf I}
\newcommand{\tildeh}{{\tilde{h}}}
\newcommand{\CH}{\text{CH}} 
\newcommand{\ch}{\textbf{ch}} 

\newcommand{\Dplane}{({\mathbb P}^2)^*}
\newcommand{\rationalmap}{\; \; - \rightarrow} 
\newcommand{\image}{\text{image} \, }
\newcommand{\pp}{\mathfrak p}
\newcommand{\mm}{\mathfrak m}
\newcommand{\Hex}{\mathcal {H}}
\newcommand{\Hexc}{\mathcal {H}^\circ}
\newcommand{\Blowup}{\text{Bl}} 
\newcommand{\spec}{\text{Spec} \,} 

\newcommand{\rank}{\text{rank}}
\renewcommand{\proof}{{\sc Proof. \;}}
\newcommand{\partltr}{{\mathfrak P}}
\newcommand{\redA}{\textcolor{red}{A}}
\newcommand{\redB}{\textcolor{red}{B}}
\newcommand{\redC}{\textcolor{red}{C}}
\newcommand{\blueD}{\textcolor{blue}{D}}
\newcommand{\blueE}{\textcolor{blue}{E}}
\newcommand{\blueF}{\textcolor{blue}{F}}
\newtheorem{Theorem}{Theorem}[section]
\newtheorem{Lemma}[Theorem]{Lemma}
\newtheorem{Proposition}[Theorem]{Proposition}

\newcommand{\arr}[6]{\left[ \begin{array}{ccc} #1 & #2 & #3\\ #4 & #5 & #6 \end{array} \right]}
\newcommand{\pasc}[6]{\left\{ \begin{array}{ccc} #1 & #2 & #3\\ #4 &
     #5 & #6 \end{array} \right\}}
\begin{document} 
\title{Degenerations of Pascal lines} 
\author{Jaydeep Chipalkatti}
\thanks{Jaydeep Chipalkatti, Department of Mathematics,
  University of Manitoba, Winnipeg, MB R3T 3R5, Canada. \\
  E-mail: \href{mailto:jaydeep.chipalkatti@umanitoba.ca}{jaydeep.chipalkatti@umanitoba.ca}}
\author{Sergio Da Silva}
\thanks{Sergio Da Silva, Department of Mathematics,
  McMaster University, Hamilton, ON L8S 4L8, Canada. \\
  E-mail: \href{mailto:dasils19@mcmaster.ca}{dasils19@mcmaster.ca}}

\maketitle

\bigskip

\parbox{16.5cm}{ \small
{\sc Abstract:} Let $\conic$ denote a nonsingular conic in the complex
projective plane. Pascal's theorem says that, given six distinct points
$A,B,C,D,E,F$ on $\conic$, the three intersection points $AE \cap BF,
AD \cap CF, BD \cap CE$ are collinear. The line containing them is
called the Pascal line of the sextuple. However, this construction may
fail when some of the six points come together. In this paper, we
find the indeterminacy locus where the Pascal line is not
well-defined and then use blow-ups along polydiagonals to
define it. We analyse the geometry of Pascals in these 
degenerate cases. Finally we offer some
remarks about the indeterminacy of other geometric elements in
Pascal's \emph{hexagrammum mysticum}. } 

\bigskip 


\medskip 


\medskip 

\tableofcontents

\section{Introduction}  
Pascal's theorem is one of the most elegant results in classical
projective geometry. Given a collection of six distinct points on a
conic, it allows us to define a highly symmetrical configuration called the
\emph{hexagrammum mysticum}. We begin with an elementary introduction to this subject; the main
results will be described in 
Section~\ref{section.results} after the required notation is available. 

\newsavebox{\pascal} 

\subsection{} \label{section.defn.pascal} 
Let $\conic$ denote a nonsingular conic in the complex projective
plane $\proj^2$. Given six distinct points $A, B, C, D, E, F$ on
$\conic$, one can arrange them into an array $\arr{A}{B}{C}{F}{E}{D}$. Then 
Pascal's theorem says that the three cross-hair intersection points 
\[ AE \cap BF, \quad AD \cap CF, \quad BD \cap CE, \] 
(corresponding to the $2 \times 2$ minors of the array) are collinear (see
Diagram~\ref{diagram.pascal.theorem}). 
\savebox{\pascal}{$\pasc{A}{B}{C}{F}{E}{D}$} 
\begin{figure}
\includegraphics[width=8cm]{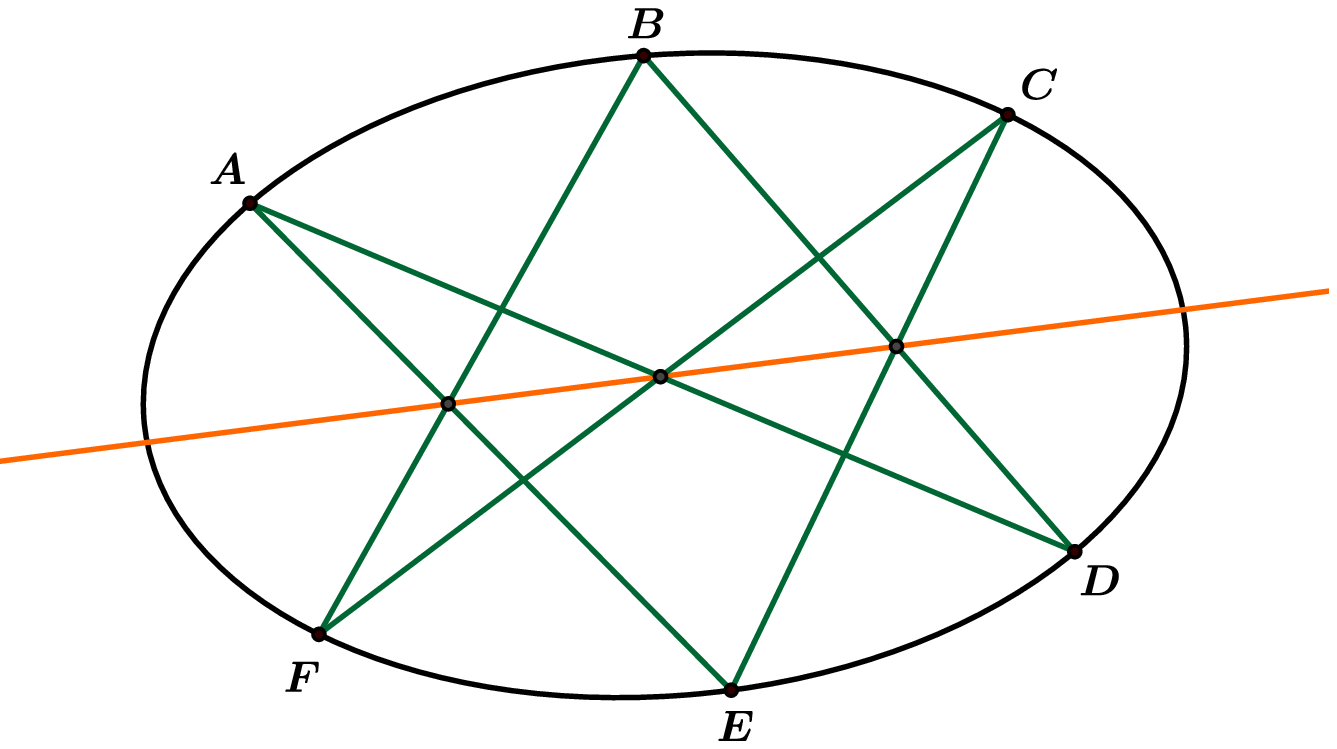}
\caption{\small The Pascal~\usebox{\pascal}}
\label{diagram.pascal.theorem} 
\end{figure} 
The line containing them is called the Pascal line (or just the
Pascal) of the array; we will denote it by
$\pasc{A}{B}{C}{F}{E}{D}$. 

\savebox{\pascal}{$\pasc{A}{B}{C}{F}{A}{D}$} 
\begin{figure}
\includegraphics[width=11cm]{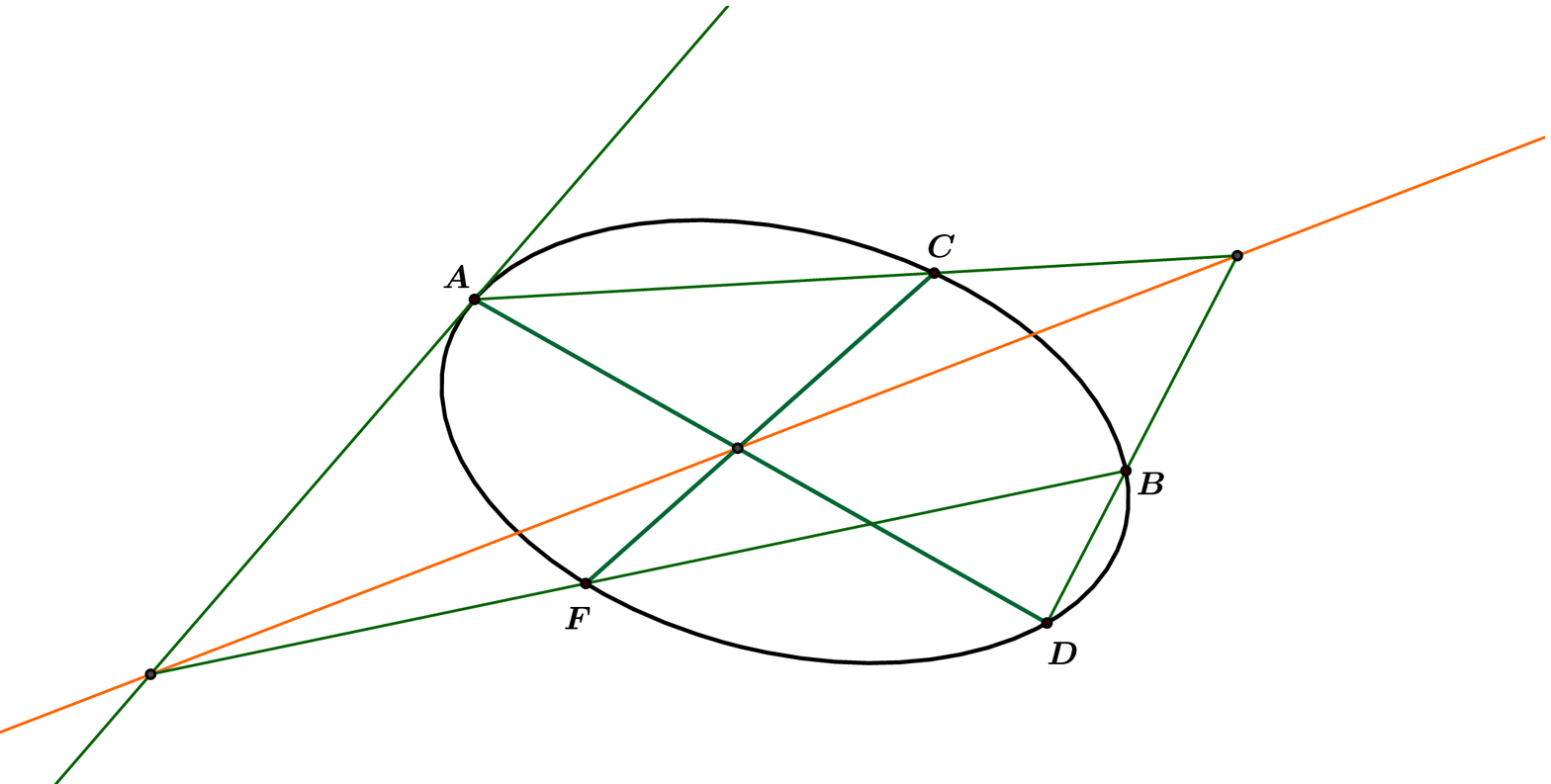}
\caption{\small The Pascal~\usebox{\pascal}} 
\label{diagram.pascal.doublepoint} 
\end{figure} 

The Pascal remains unchanged if we shuffle the rows or columns of the 
array; thus we have $12$ different ways 
\[ \pasc{A}{B}{C}{F}{E}{D} = \pasc{F}{E}{D}{A}{B}{C} =
\pasc{F}{D}{E}{A}{C}{B}  \quad \text{etc} \]
of denoting the same Pascal. Any essentially different arrangement of
the same points, such as $\pasc{D}{A}{B}{F}{C}{E}$, will generally correspond to a
different Pascal. Thus there are $6!/12 = 60$ notionally different
Pascals. It is a theorem due to Pedoe~\cite{Pedoe} that these sixty
lines are pairwise distinct if the initial six points $A, \dots, F$
are chosen in general position.

\subsection{} 
If exactly two of the points amongst $A, \dots, F$
coincide, then all the Pascals remain well-defined as
long as we follow a natural convention. If, say $P=Q$, 
then we should interpret $PQ$ as the tangent line to
the conic at $P$. (Henceforth, we denote this tangent by $\bT_P$.) For
instance, in the case of the Pascal $\pasc{A}{B}{C}{F}{E}{D}$, 
up to relabelling there are three possibilities for two of the points to coincide;
namely $A=B, A=F$ or $A=E$. In the first case, the Pascal 
 is simply the line 
$AD$, and in the second case it is the line joining $A$ to $BD \cap
CE$. The third case is shown in
Diagram~\ref{diagram.pascal.doublepoint}. 
\savebox{\pascal}{$\pasc{A}{B}{C}{F}{A}{D}$}

\subsection{}
However, things start breaking down if three points
coincide, say $A=B=C$.  Then  it is no longer obvious how to define the
Pascal $\pasc{A}{A}{A}{F}{E}{D}$, since all the three cross-hair
intersections are at $A$. Similarly, if two
pairs of points coincide, say $A=B$ and $E=F$, then it is not obvious how
to define $\pasc{A}{A}{C}{E}{E}{D}$.

\subsection{Results} \label{section.results} 
In the next section, we will introduce a projective variety $\Hex$ which acts as a parameter
space for labelled sextuples of points on the conic. Given a formal
arrangement of letters such as $s = \pasc{\bA}{\bB}{\bC}{\bF}{\bE}{\bD}$, the definition of
the Pascal line will correspond to a rational map
\[ p_s: \Hex \rationalmap \Dplane  \]
from the space of sextuples to the space of lines in the projective plane.
\begin{itemize}
  \item In Proposition~\ref{prop.indet} below, we will
characterise the indeterminacy scheme of this map. It will turn out
to be a union of \emph{polydiagonals} in $\Hex$.
\item We will remove a codimension three subvariety from $\Hex$, and blow-up
the resulting space $\Hexc$ along certain polydiagonals. The result is a
diagram
\[
  \begin{tikzcd}
    \bbX \arrow{d}{\beta} \arrow{dr}{q_s} \\ 
   \Hexc \arrow[r,dashed,"p_s"] & \Dplane 
  \end{tikzcd}
\]
where $\beta$ is a sequence blow-ups, and $q_s$ is a regular map. In
other words, $q_s$ resolves the indeterminacy in the rational map
$p_s$ (restricted to $\Hexc$).
This is the content of Theorem~\ref{theorem.resindet.Pascal}. 
\item
It is of geometric interest to see how the Pascals behave on the
exceptional loci of the blow-ups. This analysis is carried out in
Section~\ref{section.pascalmorphism.exc}. There are three cases to
consider, amongst which the partition $(2,2,2)$ is the most
symmetric and leads to the most elegant results. 
\item
 The sixty Pascals are part of a larger geometric configuration
 called the~\emph{hexagrammum mysticum}. It is comprised of $35$ more
 lines (apart from the Pascals) and $95$ more points. The questions of definability
 can also be raised for these elements. Although we do not venture
 into this analysis deeply, we offer some remarks about it in
 Section~\ref{section.HM}.
\end{itemize}

\subsection{References} 
The literature on Pascal's theorem is enormous. The standard
classical reference is by George Salmon
(see~\cite[Notes]{SalmonConics}). For later development of this
material, it is convenient to introduce the so-called `dual notation'
which makes crucial use of the outer automorphism of the symmetric
group $S_6$ (see ~\cite{Howard_etal}).
This notation, together with a host of results discovered by Cremona and
Richmond are explained by H.~F.~Baker in his note `On the 
\emph{Hexagrammum Mysticum} of Pascal' in~\cite[Note II,
pp.~219--236]{Baker}. One of the best modern surveys of this material is by Conway and
Ryba~\cite{ConwayRyba}. 
We refer the reader to~\cite{Coxeter, KK, Seidenberg} for foundational notions in 
projective geometry, and to~\cite{EisenbudHarris, Harris, Hartshorne} for those in 
algebraic geometry. In particular, we will use the notion
of a blow-up which is extensively discussed in
\cite[Ch.~IV.2]{EisenbudHarris} and \cite[Ch.~II.7]{Hartshorne}. 

\section{Partitions and polydiagonals} 
\subsection{} In this section we will introduce the necessary
geometric set-up. 
Let $\ltr$ denote the set of letters $\{\bA, \bB, \bC,
\bD, \bE, \bF \}$, and let
$\partltr$ denote the set of partitions of $\ltr$. For the sake of 
readability, we will write the partition 
\[ \pi = \{\{\bA, \bC, \bD\}, \{\bB\}, \{\bE, \bF\} \} \in \partltr \]
as $\bA\bC\bD\cdot \bB \cdot \bE\bF$, and we will say that it is of type 
$(3,2,1)$. Define $n(\pi)$ to be the cardinality of $\pi$ (which is $3$ in this case). 

Now let $\Hex = \conic^\ltr$ denote the set of maps $\ltr
\stackrel{h}{\lra} \conic$. Of course, $\Hex$ is a projective variety
isomorphic to $(\proj^1)^6$. Given $h \in \Hex$, we will usually write 
\[ h(\bA) = A, \quad h(\bB) = B, \; \text{etc} \] 
to denote the corresponding points on $\conic$.

Every $h \in \Hex$ determines a partition $\theta(h) \in \partltr$ by the
rule that $x, y \in \ltr$ belong to the same element of $\theta(h)$, 
if and only if $h(x) = h(y)$. In particular, $h$ is injective exactly
when $\theta(h) = \bA. \bB. \bC. \bD. \bE. \bF$.

\subsection{} 
For every partition $\pi \in \partltr$, we have a \emph{polydiagonal} inside $\Hex$, defined as 
follows: 
\[ \Delta_\pi: = \{ h \in \Hex: h(x) = h(y) \; \; \text{if $x$ and $y$
  belong to the same element of $\pi$} \}. \] 
Note that the defining condition says `if', and not `iff'. 
For instance, if $\pi = \bA\bC\bD\cdot \bB \cdot \bE\bF$, then $\Delta_\pi$ is the set of maps 
$\ltr \stackrel{h}{\lra} \conic$ which satisfy 
\[ h(\bA) = h(\bC) = h(\bD) \quad \text{and} \quad h(\bE) = h(\bF). \] 
Thus $\Delta_\pi$ is a nonsingular closed subvariety in $\Hex$. It is
isomorphic to $(\proj^1)^{n(\pi)}$. 

\subsection{} Partitions are partially ordered by refinement. We will
write $\pi_1 \leqslant \pi_2$ if $\pi_2$ is a refinement of $\pi_1$, in which
case $\Delta_{\pi_1} \subseteq \Delta_{\pi_2}$. In particular, the
smallest polydiagonal $\Delta_\pi$ is the one for which $\pi =
\bA\bB\bC\bD\bE\bF$. At the other end, if $\pi =
\bA. \bB. \bC. \bD. \bE. \bF$ is the trivial partition, then
$\Delta_\pi = \Hex$.

If $h \in \Delta_\pi$, then $\theta(h) \leqslant \pi$ and hence 
$\Delta_{\theta(h)} \subseteq \Delta_\pi$. For a general point $h \in
\Delta_\pi$, we have $\theta(h) = \pi$.

\subsection{} \label{section.indetcalculation.123} 
Let $[z_0,z_1,z_2]$ be the homogeneous coordinates in 
$\proj^2$. Then, for instance, the coordinates of the line $3 \, z_0 + 5 \, z_1 + 7 \,
z_2 = 0$ will be written as $\langle 3, 5, 7 \rangle$.  Identify $\conic$
with the conic $z_0 \, z_2 = z_1^2$, and fix an isomorphism 
$\tau: \proj^1 \lra \conic$ by the formula 
$\tau(a) = [1,a,a^2]$ for $a \in \complex$, and $\tau(\infty) =
[0,0,1]$. 

Choose indeterminates $a, \dots, f$, and write 
\begin{equation} A = \tau(a), \dots, F = \tau(f).
\label{points.AthroughF} \end{equation} 

Define a Pascal symbol to be an array such as 
$s = \pasc{\bE}{\bC}{\bA}{\bB}{\bF}{\bD}$, determined up to row and
column shuffles. There are sixty such symbols. 

\subsection{} \label{section.pascal.formula} 
For the moment, we fix the symbol $s =
\pasc{\bA}{\bB}{\bC}{\bF}{\bE}{\bD}$.  Given the points $A, \dots, F$ as
in (\ref{points.AthroughF}), it is easy to calculate the coordinates
of the lines $AE, BF$ etc, and eventually those
of the Pascal $\pasc{A}{B}{C}{F}{E}{D}$. They turn out to be 
$\langle u_0, u_1, u_2 \rangle$, where 
\begin{equation} \begin{aligned} 
u_0 & = abde-abdf-acde+acef+bcdf-bcef, \\ 
u_1 & = -abe+abf+acd-acf+adf-aef-bcd+bce-bde+bef+cde-cdf, \\ 
u_2 & = -ad+ae+bd-bf-ce+cf . \end{aligned} \label{coord.pascal} \end{equation} 
(This computation was programmed by us in {\sc Maple}.)

Given the polynomial ring $R = \complex[a,b,c,d,e, f]$, we have a
rational map
\[ \spec R \stackrel{p_s}{\rationalmap} \Dplane, \qquad (a, \dots, f) \rightarrow
  \langle u_0, u_1, u_2 \rangle. \] 
The indeterminacy scheme of this map is defined by the ideal $\fraka_s =
(u_0, u_1, u_2) \subseteq R$. Now it is
straightforward to calculate the minimal primary decomposition of $\fraka_s$;
in fact all of its primary components turn out to be prime
ideals. (This computation was done in {\sc Macaulay2}.) We have 
\begin{equation}  \fraka_s = \bigcap\limits_{i=1}^6 \; \pp_i,
\label{a.decomposition}   \end{equation}
where $\pp_1, \dots, \pp_6$ are the following prime ideals: 
\[ \begin{array}{lll} 
\pp_1 = (b-c,a-c), & \pp_2 = (e-f,d-f), \\ 
\pp_3 = (e-f,a-b), & \pp_4 = (d-e,b-c), & \pp_5 = (d-f,a-c), \\
\pp_6 =  (c-d,b-e,a-f). 
   \end{array} \]

 The construction above globally corresponds to a rational map
 \[ p_s: \Hex \rationalmap \Dplane. \]
 Let $\Omega_s \subseteq \Hex$ denote its indeterminacy scheme. The
preceding calculation proves the following: 
\begin{Proposition} \rm 
The scheme $\Omega_s$ is equal to the
union $\bigcup \Delta_\pi$, where $\pi$ ranges over
the following six partitions: 
\begin{equation} \label{sixpartitions.list} 
\begin{array}{lllll} 
  \vartheta_1 = \bA \bB \bC.\bD. \bE. \bF, & & 
                                             \vartheta_2 = \bA.\bB.\bC.\bD \bE \bF, \\ 
  \vartheta_3 = \bA \bB.\bC.\bD.\bE \bF, & & 
                                           \vartheta_4 = \bA.\bB
                                           \bC.\bD \bE.\bF, & & 
                                                              \vartheta_5
                                                              = \bA
                                                            \bC.\bB.\bD
                                                            \bF.\bE,
  \\ 
\vartheta_6 = \bA \bF.\bB \bE.\bC\bD. \end{array} \end{equation} 
\label{prop.indet} 
In particular, this is a reduced scheme with six irreducible
components. \end{Proposition} 

These partitions encode a simple geometric pattern, which says that 
the Pascal $\pasc{A}{B}{C}{F}{E}{D}$ can be undefined in the following three ways: 
\begin{itemize} 
\item 
all three points in either of the two rows become equal, or 
\item 
any two of the columns become equal, or 
\item 
the two rows become equal. 
\end{itemize} 
The first situation is captured by $\vartheta_1$ and $\vartheta_2$ which
are of type $(3,1,1,1)$, the second by
$\vartheta_3,\vartheta_4,\vartheta_5$ which are of type $(2,2,1,1)$,
and the last by $\vartheta_6$ which is of type $(2,2,2)$. For later
use, observe that if $i,j$ are distinct indices between $1$ and $6$, then
$\Delta_{\vartheta_i} \cap \Delta_{\vartheta_j}$ is
a polydiagonal of type $(6),(3,3), (3,2,1)$ or $(4,2)$. For instance,
\[ \Delta_{\vartheta_1} \cap \Delta_{\vartheta_3} =
  \Delta_{\bA\bB\bC.\bD.\bE\bF} \]
is of type $(3,2,1)$.

\subsection{} \label{section.indet.local} 
Assume that $h \in \Hex$ belongs to exactly one of the polydiagonals from
(\ref{sixpartitions.list}), say $\Delta_{\vartheta_1}$. Then the
corresponding maximal ideal $\mm_h$ in $R$ contains $\pp_1$, but not
$\pp_2, \dots, \pp_6$. For $2 \leqslant i \leqslant 6$, choose an
element $g_i \in \pp_i \setminus \mm_h$, and let $g = \prod
g_i$. Localising at $g$, we have an equality of ideals
$(\fraka_s)_g = (\pp_1)_g$ in $R_g$. Thus, inside the open set $\spec
R_g \subseteq \spec R$ containing $h$, the indeterminacy scheme
coincides with $(\pp_1)_g$. 

Of course, we have a similar rational map
\[ p_s: \Hex \rationalmap \Dplane, \]
for any Pascal symbol $s$. Its indeterminacy
scheme $\Omega_s$ is obtained by appropriately permuting the letters $\bA, \dots, \bF$.

\subsection{} \label{example.twosymbols}
As an example, consider the symbols
\[ s_1 = \pasc{\bA}{\bE}{\bD}{\bF}{\bB}{\bC} \quad \text{and} \quad
  s_2 = \pasc{\bA}{\bF}{\bC}{\bD}{\bE}{\bB}. \]
Now let $h$ be a general point of $\Delta_{\bA\bF. \bB\bE.\bC\bD}$;
that is to say, a map $\ltr \stackrel{h}{\longrightarrow} \conic$ such that
\[ h(\bA) = h(\bF) = A, \quad h(\bB) = h(\bE) = B, \quad h(\bC) =
  h(\bD) = C, \]
where $A, B, C$ are distinct points on the conic. Then $s_1$ leads to
the undefined Pascal $\pasc{A}{B}{C}{A}{B}{C}$ since the two rows become
equal. However, $s_2$ becomes the Pascal $\pasc{A}{A}{C}{C}{B}{B}$ which remains
well-defined (and in fact equals the line $AB$). Thus
$\Delta_{\bA\bF. \bB\bE.\bC\bD}$ is contained in $\Omega_{s_1}$, but
not in $\Omega_{s_2}$.

On the other hand, the smaller polydiagonal
$\Delta_{\bA\bC\bD\bF.\bB\bE}$ is contained in $\Omega_{s_2}$, since
its generic point leads to the undefined Pascal
$\pasc{A}{A}{A}{A}{B}{B}$ whose last two
columns are equal. 
  
\subsection{} It is a little unsatisfactory that the proof of 
Proposition~\ref{prop.indet} should rely upon a machine
computation. In this section, we show how to prove 
a weaker version of the result using only elementary algebra. It will not
be needed in the rest of the paper.

Define the expressions
\[ \begin{array}{ll} 
     P_0=bf-be+ce, &  Q_0 =(e - f) \, (b - c)  \, (ae + bd - bf - ce), \\ 
     P_1=-(c+f), & Q_1=(f-e) \, (b-c) \, (a-c+d-f), 
                   \end{array} \] 
                 and
                 \[ \delta = \det \underbrace{\left[ \begin{array}{ccc} 1 & 1 & 1
                                           \\ a & b & c \\ f & e &
                                                                   d \end{array}
                                                               \right]}_M. \]
Then it is straightforward to check that 
\[ u_0 = P_0 \, \delta + Q_0, \quad
  u_1 = P_1 \, \delta + Q_1, \quad u_2 = \delta. \]
Define the conditions
\[ \begin{array}{lll} 
  \textC_1: \; a=b=c, & \textC_2: \; d = e = f, \\
  \textC_3: \; a = b, e = f, & \textC_4: \; b=c, d=e, &\text C_5: \; a=c,
                                                 d=f, \\
     \textC_6: \; a=f, b=e, c=d. \end{array} \]
These conditions are exactly parallel to the definitions of $\pp_i$ and $\vartheta_i$
in the previous section. The following is a set-theoretic (and hence weaker) version of
Proposition~\ref{prop.indet}. 
 \begin{Proposition} \rm Given complex numbers $a, b, c, d, e, f$, we have
   $u_0 = u_1 = u_2 =0$ if and only if at least one of the conditions
   $\textC_i$ is satisfied. 
 \end{Proposition}

 \proof The `if' part follows by a
 simple computation. Any of the $\textC_i$ makes the
 matrix $M$ have rank $\leqslant 2$, forcing $\delta =0$. Furthermore,
 we have $Q_0=Q_1=0$ implying $u_i=0$ for all $i$. 
  
Now assume that $u_0=u_1=u_2=0$, and hence $Q_0=Q_1=\delta
=0$. This implies that $\rank(M) \leqslant 2$. Now
assume that none of the conditions $\textC_i$ for $1 \leqslant i \leqslant 5$ are
satisfied. Since $\textC_1$ is false, the second row of $M$ is not a multiple
of the first. Hence the third row must be a linear combination of the first
two; that is to say,
\[ f = r \, a + s, \quad e = r \, b + s, \quad d = r \, c + s, \]
for some $r,s \in \complex$. Substituting this into the $Q_i$, we get
\[ Q_0 = r \, s \, (a - b) (b - c) \, (c-a) =0 \quad \text{and} \quad
  Q_1 = r \, (r-1) \, (a - b) \, (b - c) \, (c-a) =0. \]
Now $a \neq b$,  since otherwise $f=e$ and then $\textC_3$ would
hold. A similar argument for
$\textC_4$ and $\textC_5$ shows that $a,b,c$ are pairwise distinct. Hence
$r \, s = 0$, and $r \, (r-1) =0$. If $r=0$, then $d=e=f$ which is
disallowed. Hence we must have $r=1, s=0$, which implies that
$\textC_6$ must hold. \qed

\medskip 

This implies the set-theoretic equality
\[ V(\fraka_s) = \bigcup\limits_{i=1}^6 \; V(\pp_i), \]
but we will need the stronger version in~(\ref{a.decomposition}). 

\subsection{} 
Given a Pascal symbol $s$, we should like to resolve the 
indeterminacy in the definition of the rational map
\[ p_s: \Hex \rationalmap \Dplane. \] 
That is to say, we want to construct a proper birational morphism $\bbX \longrightarrow
\Hex$ such that there is a well-defined map 
$\bbX \longrightarrow \Dplane$ which factors through $p_s$. There are possibly many ways to
do this, but we mention two which do not turn out to be
geometrically feasible.

\begin{itemize}
\item
  According to the standard formalism of \cite[Ch.~II.7]{Hartshorne},
  we can take $\bbX$ to be the blow-up of   $\Hex$ along the subscheme
  $\Omega_s$. However, since $\Omega_s$ has several components with multiple
  intersections between them, the resulting space would be highly
  singular and unwieldy to work with.
\item
  Another possibility is to construct $\bbX$ in stages, first by 
  blowing up the smallest polydiagonal $\Delta_{\bA\bB\bC\bD\bE\bF}$,
  followed by blowing up the proper transforms of the next smallest
  polydiagonals and so on. This is the central idea behind the
  configuration space constructed by Ulyanov~\cite{Ulyanov}. However,
  after trying out this approach we encountered several obstacles.
  Already at the first stage, there is a `problematic'
  locus inside the exceptional hypersurface on which the Pascal is not
  well-defined. (In a nutshell, it arises due to the appearance of
  $\delta$ in the formulae for $u_i$.) This problem recurs at 
  several intermediate stages;
  which indicates that the total number
  of blow-ups needed would be very large, and their combinatorics
  would be difficult to control.
\end{itemize}

The approach we have chosen is to remove certain polydiagonals
in $\Hex$, and then to blow-up the open sublocus $\Hexc \subseteq \Hex$
along the remaining ones. In other words, we resolve the
indeterminacies of the restricted map
\[ p_s: \Hexc \rationalmap \Dplane. \]
As mentioned above, this turns out to be geometrically the most
pragmatic solution. Our construction is uniform in the sense that it
simultaneously works for all Pascal symbols $s$.

\section{Resolution of indeterminacy}

Recall the following connection between
blow-ups and the extension of rational maps: 
\begin{Proposition} \rm 
  Let $\bV$ be a variety with a rational map $f: \bV \rationalmap 
  \proj^n$. Let $\Sigma$ be the indeterminacy scheme (also known as the 
  scheme of base points) of $f$. Let $\Blowup_\Sigma(\bV)$ denote the
  blow-up of $\bV$ along $\Sigma$. Then $f$ extends to a regular 
  morphism 
  \[ \widetilde{f}: \Blowup_\Sigma(\bV) \longrightarrow \proj^n. \]
\label{prop.ext.rational} \end{Proposition}

\proof See \cite[Ch.~2.II, Example~7.17.3]{Hartshorne}. \qed 

This result will be used in the proof of
Theorem~\ref{theorem.resindet.Pascal} below. In summary, given an
affine chart $\spec R \subseteq \bV$, the map $f$ is given by a projective
$(n+1)$-tuple of functions $[v_0, \dots, v_n]$ such that the $v_i$
generate the ideal of $\Sigma \cap \spec R$. By construction, the pullback of $\Sigma$ to
$\Blowup_\Sigma(\bV)$ is Cartier (i.e., locally defined by
a single nonzero divisor). This divisor can be `cancelled out' from
the tuple, which gives a well-defined expression for
$\widetilde{f}$. The reader will find many such examples
in~\cite[Ch.~IV.2]{EisenbudHarris}. A thematically similar discussion
is given in~\cite[Ch.~2.II]{Hartshorne} preceding the proof of the
proposition. 

\subsection{} We now proceed in the following steps: 
\begin{itemize}
\item
  Define the open set 
\[ \Hexc  = \Hex  - \bigcup\limits_{\mu}
  \; \Delta_\mu,  \]
where the union is over all polydiagonals of type $(3,2,1)$ or
$(4,1,1)$. Notice that the union automatically includes all polydiagonals of type $(5,1),
(3,3), (4,2)$ and $(6)$.
\item
If $\pi$ is of type $(2,2,1,1), (3,1,1,1)$ or $(2,2,2)$, then define
$\Delta_\pi^\circ = \Hexc \cap \Delta_\pi$, which we call an open
polydiagonal of type $\pi$. 
\item
Now let $Z \subseteq \Hexc$ denote the union of all open 
  polydiagonals of type $(2,2,2)$, and
 let $Y \subseteq \Hexc$ denote the union of all open 
 polydiagonals of type $(2,2,1,1)$ or $(3,1,1,1)$. Then we have
 \[ Z \subseteq   Y \subseteq \Hexc. \]
\end{itemize}

\begin{figure}
\includegraphics[width=14cm]{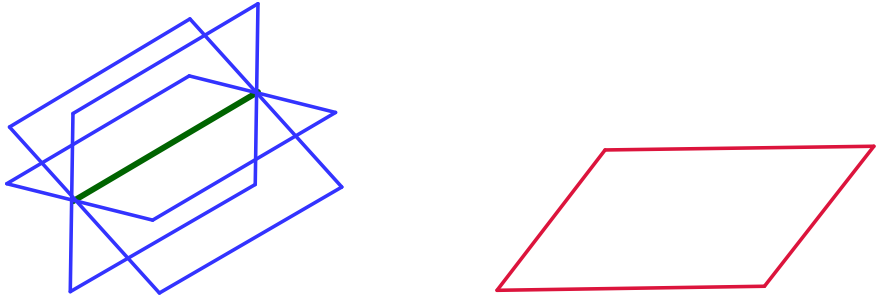}
\caption{Polydiagonals in $\Hexc$. The green line has
  type $(2,2,2)$ and each of the three blue rectangles has type
  $(2,2,1,1)$. The red rectangle has type $(3,1,1,1)$.}
\label{diagram.blowup}
\end{figure}

A schematic picture is shown in Diagram~\ref{diagram.blowup}. The
thick green line represents a typical open polydiagonal of type
$(2,2,2)$, such as $\Delta_{\bA\bF.\bB\bE.\bC\bD}^\circ$.
It is the transverse intersection of three open polydiagonals of
type $(2,2,1,1)$, shown as blue rectangles. In this example, 
\begin{equation}
  \Delta_{\bA\bF.\bB\bE.\bC\bD}^\circ =
  \Delta_{\bA\bF.\bB\bE.\bC.\bD}^\circ
  \cap \Delta_{\bA.\bF.\bB\bE.\bC\bD}^\circ 
  \cap \Delta_{\bA\bF.\bB.\bE.\bC\bD}^\circ. \label{intersection.2211}
\end{equation} 
  The red rectangle represents a typical open polydiagonal of type
  $(3,1,1,1)$. Thus $Z$ is the union of all green lines, and $Y$
  is the union of all blue and red rectangles.

  In general, a point $h \in Y$ will be in the indeterminacy locus of
  $p_s$ for some Pascal symbols $s$ and not for others. The issue is decided
  by how the combinatorial structure of $\theta(h)$ interacts with
  that of $s$ (cf.~Section~\ref{example.twosymbols}). 

  \begin{Lemma} \rm
    Let $X'$ be the blow-up of $\Hexc$ along
    $\Delta_{\bA\bF.\bB\bE.\bC\bD}^\circ$. Then the proper transforms of the
    three polydiagonals appearing on the right-hand side
    of~(\ref{intersection.2211}) are pairwise disjoint in $X'$. 
  \label{lemma.3rectangles} \end{Lemma}
  \proof
The result will follow from a local calculation. We should like to blow up the
ideal $J = (f-a,e-b,d-c)$ inside the ring $R =
\complex[a,b,c,d,e,f]$. The resulting space is covered by three affine charts
corresponding to the generators of $J$. For instance, choose a parameter $w = f-a$, and
let $e-b = q_1 \, w, \, d-c = q_2 \, w$. The blow-up locally
corresponds to the ring map
\[ R \longrightarrow \underbrace{\complex[a,b,c,w, q_1, q_2]}_S, \]
which sends $a,b,c$ to themselves, together with
\[ f \rightarrow a+ w, \quad 
  e \rightarrow b+ q_1 \, w, \quad 
  d \rightarrow c + q_2 \, w. \] 
The three polydiagonals respectively correspond to the ideals
$(f-a, e-b), (e-b,d-c)$ and $(f-a, d-c)$ in $R$. Since their proper transforms are
the $S$-ideals $(1), (q_1, q_2)$ and $(1)$ respectively, only the
middle one has nonempty support in the affine chart
$\spec S$. The calculation for the other two charts is essentially
identical, which proves the lemma. \qed 

\medskip 

Pictorially, the lemma says that when the green line is blown up, the
proper transforms of the three blue rectangles get separated. 

\subsection{} \label{section.markedpoints}
  Now let $\bbX' = \Blowup_Z(\Hexc)$ denote the blow-up of $\Hexc$ along
  $Z$, with structure morphism
  \[ \bbX' \stackrel{\beta_1}{\longrightarrow} \Hexc. \]
  Let $Y' \subseteq \bbX'$ denote the proper transform of
  $Y$. Notice that, before the blow-up, any of the red rectangles was already disjoint from
  any of the blue ones, since we have removed their intersections from
  $\Hex$ while arriving at $\Hexc$. Hence $Y'$  is a disjoint union of
  varieties which are individually easy to handle. 

If $h$ is a point in $\Delta_{\bA\bF.\bB\bE.\bC\bD}$, then
  $\beta_1^{-1}(h)$ is isomorphic to $\proj^2$. This plane
intersects the proper transform of each of the blue rectangles in a
point. In order to represent these three points, it will be convenient to
have a more symmetric affine chart inside the blow-up. In the notation
of Lemma~\ref{lemma.3rectangles}, choose the parameter $t = (f+e+d)-(a+b+c)$ and
let $f-a = p_1 \, t$ and $e-b = p_2 \, t$. The blow-up locally
corresponds to the ring map
\[ R \longrightarrow \complex[a,b,c,t, p_1, p_2], \]
which sends $a,b,c$ to themselves, together with
\[ f \rightarrow a+ p_1 \, t, \quad 
  e \rightarrow b+ p_2 \, t, \quad 
  d \rightarrow c + t - p_1 \, t - p_2 \, t. \] 
The proper transforms of the three polydiagonals correspond to the ideals $(p_1, p_2), (p_2, 1-p_1-p_2)$
and $(p_1, 1-p_1-p_2)$. Hence the three marked points are given by
\begin{equation}
  W_{\bA\bF.\bB\bE}: p_1=p_2 =0, \qquad
  W_{\bB\bE.\bC\bD}: p_1=1, \, p_2 =0, \qquad
  W_{\bA\bF.\bC\bD}: p_1=0, \, p_2 = 1. 
\label{3marked.points} \end{equation} 
Of course, similar statements hold for any partition of type $(2,2,2)$. 

\medskip

Finally, let $\bbX = \Blowup_{Y'}(\bbX')$ denote the blow-up along
  $Y'$, with structure morphism
  \[ \bbX \stackrel{\beta_2}{\longrightarrow} \bbX'. \]
  Now let $\beta = \beta_2 \circ \beta_1$, so that we have a morphism
  \[ \bbX \stackrel{\beta}{\longrightarrow} \Hexc. \] 

The next theorem says that the rational map $p_s$ extends to a
morphism $q_s$ on $\bbX$. 

\begin{Theorem} \rm Let $s$ be any Pascal symbol. Then we have a
  commutative diagram
  \[
  \begin{tikzcd}
    \bbX \arrow{d}{\beta} \arrow{dr}{q_s} \\ 
   \Hexc \arrow[r,dashed,"p_s"] & \Dplane 
  \end{tikzcd}
\]
\label{theorem.resindet.Pascal} \end{Theorem}
The commutativity has the following meaning: if $\tildeh \in \bbX$ is a
point such that $p_s$ is already defined at $h = \beta(\tildeh)$, then
$q_s(\tildeh) = p_s(h)$.

\proof 
Let $\tildeh \in \bbX$, and write $h  = \beta(\tildeh)$. If $h \in \Hexc
\setminus Y$, then $\beta$ is an isomorphism in an open neighbourhood of $\tildeh$
and hence $q_s$ is well-defined at $\tildeh$.

Now assume that $h \in Y$. If $p_s$ is already defined at $h$, then we
set $q_s(\tildeh) = p_s(h)$. Hence we may further assume that $p_s$ is
undefined at $h$. If $h \notin Z$, then in a sufficiently small open neighbourhood $U$ of
$h$ we have $\beta^{-1}(U) \simeq \Blowup_{(Y \cap U)}(U)$. By the
argument of Section~\ref{section.indet.local}, the 
indeterminacy scheme of $p_s$ coincides with $Y \cap U$ in $U$. 
Then Proposition~\ref{prop.ext.rational} implies that $p_s$ extends to a regular map over
$\beta^{-1}(U)$.

Now assume $h \in Z$. Then the same argument implies
that $p_s$ extends to a regular morphism $q_s'$ on 
$\beta_1^{-1}(U)$ for an open neighbourhood $U$ of $h$. Now define $q_s(\tildeh) =
q_s'(\beta_2(\tildeh))$. This proves the theorem. \qed

\subsection{} 
\label{section.example.pascals}
We can calculate the hitherto undefined Pascals using the recipe 
above. For instance, let $h \in \Hexc$ be such that 
\[ h(\bA) = h(\bB) = h(\bC) = \tau(3), \quad h(\bD) = \tau(1),
  \quad h(\bE) = \tau(7), \quad h(\bF) = \tau(4),  \]
in the notation of Section~\ref{section.indetcalculation.123}. 
Then the Pascal is undefined for $s =
\pasc{\bA}{\bB}{\bC}{\bF}{\bE}{\bD}$, 
since all points in the top row become equal. 

In order to define the Pascal over the fibre $\beta^{-1}(h)$, we need to blow up the ideal $I = (b-a,c-a)$ inside
the ring $R = \complex[a,b,c,e,d,f]$. Choose a parameter
$b - a = t$ and let $c - a = p \, t$. 
Locally above the point $h$, the morphism $\beta$ is given by the ring map
\[
  R \longrightarrow
  \underbrace{\complex[a,d,e,f,p,t]}_S \]
which acts by sending $a,d,e,f$ to themselves, together with 
$b \rightarrow a+t, c \rightarrow a + p \, t$. 
Now make the substitutions
\[ a \rightarrow 3, \quad
  b \rightarrow 3 + t, \quad
  c \rightarrow 3 + p \, t, \quad
  d \rightarrow 1, \quad e \rightarrow 7, \quad f \rightarrow 4,   \]
into the expressions $u_i$ from Section~\ref{section.pascal.formula}; then we get
\[ 
    u_0 = -3 \, t \, (8 \, p \, t + 3 \, p + 21), \quad 
    u_1 = 6 \, t \, (p \, t + 2 \, p + 5), \quad 
    u_2 = -3 \, t \, (p + 1).
  \]
  Since $ \langle u_0, u_1, u_2 \rangle$ is a projective triple, we may cancel the
  $t$. After substituting $t=0$ (which corresponds to the exceptional
  locus of the blow-up), we get the formula 
  \[  \langle u_0, u_1,u_2 \rangle  = \langle \, -3 \, ( 3 \, p+21), 6 \, (2 \, p+5), -3 \,
    (p+1) \, \rangle = \langle 3p+21, -4p-10, p+1 \, \rangle \]
  for the line coordinates of the Pascal, as $p$ describes a variable point in $\beta^{-1}(h)
  \simeq \proj^1$. (It is understood that the point $p =
  \infty$ will be captured in a different affine chart.) Notice the identity $ u_0 + 3 \, u_1 + 3^2 \,
  u_2=0$, which implies that the Pascal always passes through the
  triple point $\tau(3)$. In the next section, we will carry out this
  analysis more generally. 

\section{The Pascal morphism on the exceptional loci} 
\label{section.pascalmorphism.exc} 

If $h \in \Hexc$ is a sextuple which lies in one of the polydiagonals
which have been blown up, then it is of interest to see how the
Pascals behave on the fibre
\[ \beta^{-1}(h) = X_h. \]

\subsection{} \label{section.threecases}
In this section we will describe the maps
\[ q_s: X_h \longrightarrow \Dplane, \]
as $s$ runs through all Pascal symbols. There are three cases to
consider; namely when $\theta(h)$ is of type $(3,1,1,1), (2,2,1,1)$
or $(2,2,2)$. The last case is the most symmetric, and as such
geometrically the most elegant.

The first two of these partitions are of codimension two in
$\Hexc$, and $X_h \simeq \proj^1$ in these
cases. In the $(2,2,2)$ case, $X_h$ is isomorphic to $\proj^2$ blown up in
the three marked points described in Section~\ref{section.markedpoints}.
The exceptional locus of this blow-up consists
of three disjoint copies of $\proj^1$. We label them as
$L_{\bA\bF.\bB\bE},  L_{\bB\bE.\bC\bD}$ and $L_{\bA\bF.\bC\bD}$ in the
context of the example given there, and similarly for other
partitions. 

\subsection{The $(3,1,1,1)$ case}

Let $\theta(h) = \bA\bB\bC.\bD.\bE.\bF$. Fix distinct fix
points $M, P,Q,R$ on the conic, and assume that $h$ acts as follows:
\[ \bA, \bB, \bC \rightarrow M, \qquad
  \bD \rightarrow P, \qquad
  \bE \rightarrow Q, \qquad
  \bF \rightarrow R. \]

Let $s$ be a Pascal symbol. First assume that $\bA, \bB, \bC$ are not in the same row
of $s$. Then, up to the permutation of these three letters, $s$ has 
either of the two forms:
\[ \pasc{\bA}{\bB}{x}{\bC}{y}{z} \quad \text{or} \quad 
  \pasc{\bA}{\bB}{x}{y}{z}{\bC}, \]
where $\{x,y,z\} = \{\bD, \bE, \bF\}$. In either case, it follows from the 
definition of the Pascal that it is already defined at $h$. For the
first class of symbols, it is the line $Mz$ (where $z$ stands for $P,
Q$ or $R$ as the
case may be). For the second, it is the tangent line $\bT_M$ to
the conic at $M$. 

Now assume that $s = \pasc{\bA}{\bB}{\bC}{x}{y}{z}$, where $\{x,y,z\} =
\{\bD, \bE, \bF\}$. Then $h$ lies in the indeterminacy locus of
$p_s$. 
\begin{Proposition} \rm With notation as above, $q_s$ induces an
  isomorphism of the fibre $X_h$ with the pencil of lines through
  $M$.
\end{Proposition}

\proof This will follow by an explicit calculation as in
Section~\ref{section.example.pascals}. 
Assume $s = \pasc{\bA}{\bB}{\bC}{\bF}{\bE}{\bD}$. Using an
automorphism of $\proj^1$, we may assume that $M,P,Q,R$ respectively correspond to
the points $m,-1,0,1$ on $\proj^1$. Now substitute
\[ \begin{array}{lll} 
    a \rightarrow m, & b \rightarrow m + t, & c \rightarrow m+ p \,\ t,
     \\ 
     d \rightarrow -1, & e \rightarrow 0, & f \rightarrow 1, 
   \end{array} \]
into the formulae for $u_i$, factor out $t$ from the
projective triple $\langle u_0,u_1,u_2 \rangle$ and
substitute $t=0$. Then the formula for the Pascal turns out to be
\[ \lambda_p = \langle -m \, p, 2 \, m -p \, m +p, p-2 \rangle. \]
Its dot product with the vector $\tau(m) = [1,m,m^2]$ is zero, which
implies that the line passes through $M$. Since its coordinates are
linear in the parameter $p$, we get the desired isomorphism. The
other cases of $s$ follow by symmetry. \qed 

\subsection{The $(2,2,1,1)$ case}

Now let $\theta(h) = \bA\bB.\bE\bF.\bC.\bD$. Fix distinct fix
points $M, N,P,Q$ on the conic, and assume that $h$ acts as follows:
\[ \bA, \bB \rightarrow M, \qquad
  \bE, \bF \rightarrow N, \qquad
  \bC \rightarrow P, \qquad
  \bD \rightarrow Q. \]

There are four Pascal symbols $s$ for which $h$ lies in the
indeterminacy locus, namely
\[
  \pasc{\bA}{\bB}{x}{\bE}{\bF}{y} \quad \text{and} \quad 
  \pasc{\bA}{\bB}{x}{\bF}{\bE}{y}, \]
where $\{x,y\} = \{\bC, \bD\}$. We will give the result for one of
these patterns, and the others will follow by symmetry.

\begin{Proposition}
  \rm
  Assume $s =
  \pasc{\bA}{\bB}{\bC}{\bF}{\bE}{\bD}$. Then $q_s$ induces an
  isomorphism of the fibre $X_h$ with the pencil of lines through
  the point $MQ \cap NP$. 
\end{Proposition}

\proof
As before, this follows by an explicit calculation. We may assume that
$M, N, P, Q$ respectively correspond to $m,n,1,-1$ on $\proj^1$. The
recipe involves blowing up the ideal $(b-a, f-e)$ inside the ring $R =
\complex[a,b,c,d,e,f]$. Hence, make substitutions $b = a + t, f = e +
p \, t$ into the $u_i$ and calculate the
Pascal as above. Its line coordinates turn out to be
\[\langle m^2p - mp - n^2 - n, m^2p - 2mp + n^2 + 2n + p + 1, -mp - n
  + p - 1 \rangle,  \]
whose dot product with
$MQ \cap NP =
[n-m+2, m+n, 2 m n + m - n]$ is zero. Then the linearity in
parameter $p$ establishes the isomorphism. \qed

\smallskip

Now let us consider those $s$ for which the Pascal is already defined
at $h$. The sets $\{\bA, \bB\}, \{\bE, \bF\}$ are in symmetric
positions, and so are the letters within each set. Hence,
up to these shuffles, the essentially distinct cases are as follows:
\[
  \pasc{\bA}{\bB}{x}{y}{\bE}{\bF}, \quad 
  \pasc{\bA}{\bE}{\bF}{\bB}{x}{y}, \quad 
  \pasc{\bA}{\bE}{x}{\bB}{y}{\bF}, \quad 
  \pasc{\bA}{\bE}{x}{\bF}{\bB}{y}, \quad 
  \pasc{\bA}{\bE}{x}{y}{\bB}{\bF}, \]
where $\{x,y\} = \{ \bC, \bD\}$. The corresponding Pascals can be
found by tracking the cross-hair intersections. In the first three
cases it is the line $MN$. 
In the fourth case, it is the line joining 
$My \cap Nx$ and $Ny \cap Mx$,  where $x, y$ stand for
either $P$ or $Q$ depending on the choice of the bijection 
$\{x,y\} \rightarrow \{ \bC, \bD\}$.  In the fifth case, it is the line
joining $\bT_M \cap yN$ and $\bT_N \cap xM$. This completes the
discussion of the $(2,2,1,1)$ case.

\subsection{The $(2,2,2)$ case} \label{section.222case} 
Let $\theta(h) = \bA \bF.\bB \bE.\bC \bD$. Fix three distinct points
$P, Q, R$ on the conic, and assume that $h$ acts as follows: 
\[ \bA, \bF \lra P, \qquad \bB, \bE \lra Q, \qquad \bC, \bD \lra R. \] 
We should like to describe the morphisms $q_s: X_h \lra
\Dplane$ for varying symbols $s$.  

\begin{figure}
\includegraphics[width=14cm]{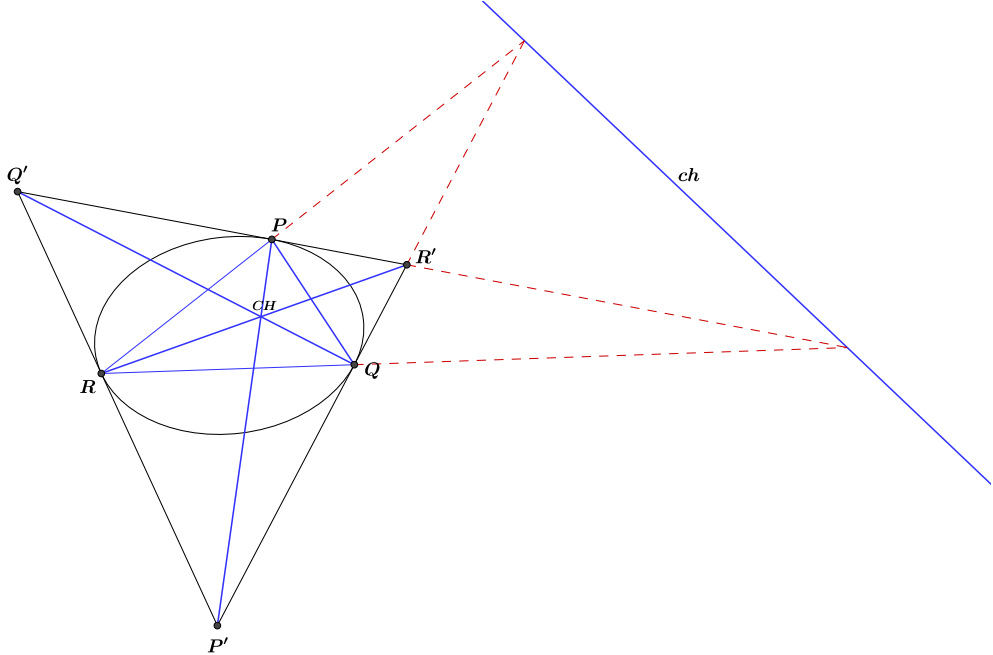}
\caption{\small Geometric elements in the $(2,2,2)$ case} 
\label{diagram.222case} 
\end{figure}

Recall that each point in the projective plane has a polar line with
respect to the conic $\conic$, and conversely each line has a pole
(see~\cite[Ch.~VI]{Seidenberg}). In Diagram~\ref{diagram.222case}, the
points $P',Q',R'$ are respectively the poles of lines $QR, PR, PQ$. Then $P'Q'R'$ is called
the polar triangle of $PQR$. It is a theorem due to Chasles that these two
triangles are in perspective (see~\cite[\S 99]{SalmonConics}); that
is to say, 
\begin{itemize} 
\item 
lines $PP', QQ'$ and $RR'$ are concurrent in a point $\CH$, and 
\item 
the points $PQ \cap P'Q', PR \cap P'R'$ and $QR \cap Q'R'$ are collinear on
a line $\ch$.
\end{itemize} 
(The intersection $PQ \cap P'Q'$ is not shown in the diagram.) Now the
behaviour of $q_s$ (restricted to fibre $X_h$) is described in the following theorem. 

\begin{Theorem} \rm 
\begin{enumerate} 
\item 
The map $q_s$ is constant for $44$ values of the symbol $s$. Amongst these, there
are $8$ values each for which the image is the line $PQ, PR$ or $QR$. There
are $4$ values each for which the image is $PP', QQ'$ or
$RR'$. Finally, there are $8$ values of $s$ for which the image is
$\ch$. 
\item 
The map $q_s$ is non-constant for $16$ values of $s$. There
  are $4$ values for which the image is the pencil of lines through
  $P$, and similarly $4$ each for $Q$ and $R$. For the remaining $4$
  values, the image is all of $\Dplane$. 
\end{enumerate} 
\label{theorem.222} \end{Theorem} 

In order to explain the patterns precisely, let the letters $x,y,z$ 
stand for elements of the three subsets $\{\bA, \bF \}, \{\bB, \bE \}, \{\bC, \bD
\}$ in some order. Moreover, if for instance $x$ stands for one of the letters $\{\bB, \bE\}$, then it
will also serve as a placeholder for the corresponding point on the conic, namely $Q$. 

For the moment, let us allot the letters arbitrarily in the following way: 
\[ x \ra \{\bA, \bF\}, \quad y \ra \{\bC, \bD\}, \quad z \ra \{\bB,
\bE\}. \] 

Now $q_s$ is a constant map in these cases: 
\begin{itemize} 
\item 
If $s = \pasc{x}{x}{y}{y}{z}{z}$, then $\image(q_s) = xz$. 
For instance, $s = \pasc{\bF}{\bA}{\bD}{\bC}{\bB}{\bE}$ obeys this
pattern, and hence the map $q_s$ is constant with image $PQ$. 
\item 
If $s = \pasc{x}{y}{z}{x}{z}{y}$, then the image is the line
$xx'$. For instance, if $s=\pasc{\bA}{\bC}{\bE}{\bF}{\bB}{\bD}$ then
it is the line $PP'$. 
\item 
For $s = \pasc{x}{y}{z}{y}{z}{x}$, it is the line $\ch$. 
\end{itemize} 

As for the variable lines, 
\begin{itemize} 
\item 
For $s = \pasc{x}{y}{y}{x}{z}{z}$, the image of $q_s$ is the pencil of lines through 
$x$. 
\item 
For $s = \pasc{x}{y}{z}{x}{y}{z}$, the image is all of $\Dplane$. 
\end{itemize} 
The same is true of any of the six possible allotments. This covers
all the sixty Pascal symbols. 

\proof
All the proofs follow by straightforward computations. We will prove three
cases for the sake of illustration.

For example, suppose that $s = \pasc{\bA}{\bB}{\bC}{\bF}{\bE}{\bD}$ which fits
the pattern $\pasc{x}{y}{z}{x}{y}{z}$. Using an automorphism of
$\proj^1$, we may assume that $P=\tau(1), Q=\tau(0), R=\tau(-1)$. Following the
generators of the ideal $I = (f-a,e-b,d-c)$ in $R = \complex[a,b,c,d,e,f]$, make
substitutions
\[ f \rightarrow a+t = 1+t, \quad
  e \rightarrow b+p_1 \, t = p_1 \, t, \quad
  d \rightarrow c + p_2 \, t = -1 + p_2 \, t \]
into the formulae for the $u_i$. Cancel the factor $t$ and
substitute $t=0$. After substituting $p_1 \rightarrow p_1/p_0, \; 
p_2 \rightarrow p_2/p_0$ and homogenizing, we see that the Pascal morphism 
is given by the composite
\[ X_h \longrightarrow \proj^2 \longrightarrow \Dplane, \]
where the first map blows down the three $L$-lines (see
Section~\ref{section.threecases}),
and the second map is the isomorphism 
\[ \left[ \begin{array}{ccc} p_0 & p_1 & p_2 \end{array} \right] \rightarrow
  \left[ \begin{array}{ccc} p_0  & p_1 & p_2 \end{array} \right]
  \;
  \left[ \begin{array}{rrr} 0 & -1 & -1 \\ -2 & 0 & 2 \\ 0 & 1 & -1
                                                                 \end{array}
                                                               \right]. \]
                                                             
As a second example, suppose that
$s = \pasc{\bA}{\bC}{\bD}{\bF}{\bB}{\bE}$ which fits the pattern
$\pasc{x}{y}{y}{x}{z}{z}$. The 
Pascal is not defined at $h$, since the second and third column
are equal. According to the recipe, we first need to blow up
the open polydiagonal $\Delta^\circ_{\bA\bF.\bB\bE.\bC\bD}$, and then the proper transform
of $\Delta^\circ_{\bA.\bF.\bB\bE.\bC\bD}$. We will follow the
notation of Section~\ref{section.markedpoints}, where the first blow-up is
already done. Since the proper transform corresponds to the ideal $(p_2, 1- p_1-p_2)$ in
$S$, we let $1-p_1 - p_2 = p_2 \, r$. Then the second blow-up corresponds
to the map
\[ \complex[a,b,c,t, p_1, p_2] \longrightarrow
  \complex[a,b,c,t,p_2, r], \]
where $a,b,c,t,p_2$ map to themselves and $p_1 \rightarrow 1-p_2 \, r
-p_2$. Now make these substitutions in the $u_i$, cancel the factor
$t \, p_2$ and set $t=0$. Using the same $P, Q, R$ as above, we get
\[ \langle u_0, u_1, u_2 \rangle = \langle -2, r, 2- r \rangle. \]
Its dot product with $\tau(1) = [1,1,1]$ is zero, hence this
represents a pencil of lines through $P$.

\medskip 

As a third example, suppose that $s =
\pasc{\bA}{\bC}{\bB}{\bD}{\bE}{\bF}$, which
fits the pattern $\pasc{x}{y}{z}{y}{z}{x}$. Then the Pascal is already
defined at $h$. By definition, it is the line passing through
\[ \bT_P \cap QR, \quad \bT_Q \cap PR \quad \text{and} \quad \bT_R \cap PQ, \]
which is $\ch$.

The other cases follow by similar calculations.
We leave the details to the industrious reader. \qed 

\section{The Hexagrammum Mysticum}
\label{section.HM}

We have seen that six points on a conic lead to a collection of
sixty Pascal lines in the plane. These lines satisfy some incidence
relations, leading to several other geometric elements such as the
Kirkman and Steiner points, and the Cayley and Salmon lines. This
entire structure is sometimes called the \emph{Hexagrammum
  Mysticum}; it is described in detail in the article by
Conway and Ryba~\cite{ConwayRyba}. The question of definability which we have
considered for Pascals can also be raised for these geometric
elements. We offer some brief remarks in this direction.
\subsection{Kirkman points}
It is a theorem due to Kirkman that the Pascals 
\begin{equation} 
  \pasc{\redA}{\blueE}{\redC}{\blueD}{\redB}{\blueF},
  \quad
  \pasc{\redB}{\blueD}{\redA}{\blueF}{\redC}{\blueE},
  \quad
  \pasc{\redC}{\blueF}{\redB}{\blueE}{\redA}{\blueD}
\label{three.pascals} \end{equation} 
are concurrent; their common point is called a Kirkman point. This
pattern is obtained by starting from the array
$\left[ \begin{array}{ccc}
   \redA & \redB & \redC \\
    \blueF & \blueE & \blueD \end{array} \right]$, and arranging the
red (respectively blue) letters in a $\vee$ shape (respectively a $\wedge$
shape) as shown. As we scan the arrays in (\ref{three.pascals}) from left to right, the red 
letters undergo a cyclic shift $ABC-BCA-CAB$, whereas the blue letters
undergo an anti-cyclic shift $DEF-FDE-EFD$. 

By shuffling the letters $\{A, \dots, F\}$ in all possible ways, 
we get altogether sixty Kirkman points in the plane. As in the case of
the Pascal line (cf.~Section~\ref{section.pascal.formula}), it
is straightforward to calculate the coordinates of the Kirkman point
described above and find the indeterminacy locus. We have carried out this computation
in {\sc Macaulay2}. It turns out that this locus is
scheme-theoretically defined by the ideal $\; \bigcap \, J_w$, where $w$ ranges over the symbols 
\[ \begin{aligned} {} & def, \; cef, \; ab.ef, \; bdf, \; ac.df, \;
    bc.de, \; ade, \; bcf, \; bd.cf, \; ae.cf, \; ad.cf, \; ce.bf, \;
    ce.bd, \; ace, \\ & ad.ce, \; ae.bf, \; ad.bf, \; ae.bd, \; abd,
    \; abc. \end{aligned} \]

This is to be interpreted as follows: if $w$ is a letter triple such as $def$, then $J_w$ stands for 
the ideal $(d-e,d-f)$. If $w$ is a dot-separated pair such as $ab.ef$, then 
$J_w$ stands for $(a-b,e-f)$. Hence the indeterminacy scheme is a union of
polydiagonals of types $(3,1,1,1)$ and $(2,2,1,1)$. In particular, all
the Kirkman points remain well-defined
as long as the initial six points are pairwise distinct or if
there is a single double point. 

\subsection{Steiner points} The Pascal lines satisfy another incidence theorem due to
Steiner. The lines 
\[ \pasc{A}{B}{C}{F}{E}{D}, 
\quad \pasc{A}{B}{C}{D}{F}{E}, 
\quad \pasc{A}{B}{C}{E}{D}{F}, \] 
are also concurrent, and their common point is called a Steiner
point. In this case, the top row is fixed at $(A,B,C)$ and the bottom
row goes through cyclic permutations of $(F,E,D)$.
We get altogether $20$ Steiner points by shuffling the labels $\{A, \dots,
F\}$. 

Now two of these Pascals may become simultaneously undefined on
certain polydiagonals, and then the corresponding Steiner point is
also undefined. However, in contradistinction to Pascal lines and Kirkman
points, a Steiner point may become undefined even if $A, \dots, F$
are pairwise distinct. In summary, the situation is as follows (see~\cite{Chipalkatti}). We will say
that a sextuple of points $A, \dots, F$ on $\conic \simeq \proj^1$ is
tri-symmetric\footnote{The rationale behind
    this term is explained in \cite{Chipalkatti}.}, if it is projectively 
equivalent to the set 
\[ \left\{0, 1, \infty, \alpha, \frac{\alpha-1}{\alpha}, \frac{1}{1-\alpha} \right\} \] 
for some complex number $\alpha$. In that case, at least one of the
Steiner points becomes undefined. 

\subsection{Cayley lines, Pl{\"u}cker lines and Salmon points}
The Kirkman points satisfy a collinearity theorem, leading to the
so-called $20$ Cayley lines. Similarly, incidences of Steiner points lead to
$15$ Pl{\"u}cker lines, and those of Cayley lines lead to
$15$ Salmon points. We will not describe the combinatorics of these incidences 
here since this is best done via the `dual notation'
(see~\cite{Baker, ConwayRyba}).

One can carry out a similar analysis for these geometric elements, and
deduce the following (see~\cite{Chipalkatti}): 
\begin{itemize} 
\item 
If the sextuple is tri-symmetric, then at least one of the Cayley lines becomes 
undefined. 
\item 
If it is tri-symmetric with $\alpha = \sqrt{-1}$, then at least one of the
Pl{\"u}cker lines and one of the Salmon points becomes undefined. 
\end{itemize} 

Thus, the indeterminacy loci for Cayley and Pl{\"u}cker lines, as well
as those for Steiner and Salmon points are intricate subvarieties
of $\Hex$ which are not confined to the polydiagonals. One would need
to carry out a detailed analysis of their geometry in order to
resolve the indeterminacies. We hope to take this up in
a possible sequel to this paper.

\medskip

{\sc Acknowledgement:} The second author was supported in part by a
PIMS postdoctoral fellowship and an NSERC postdoctoral fellowship.

\medskip 

\centerline{--} 

\end{document}